\documentclass[preprint,12pt]{elsarticle}
\usepackage{mathrsfs}
\usepackage{amsfonts}
\usepackage{amssymb,bm}
\usepackage{amsmath}
\usepackage{amsthm}
\usepackage{bbding}
\allowdisplaybreaks
\newtheorem{thm}{Theorem}[section]

\newtheorem{lem}[thm]{Lemma}

\newtheorem{cor}[thm]{Corollary}

\newtheorem{rmk}[thm]{Remark}

\newproof{pf}{Proof}
\newcommand{\qedd}{\hspace*{\fill}$\Box$\medskip}   %end box of a proof
\def\ker{\mbox{\rm{ker}}}
\def\tr{{\rm{tr}}}
\def\nor{{\rm{N}}}

\journal{}
\begin{document}

\begin{frontmatter}

\title{Complete permutation polynomials induced
from complete permutations of subfields}
\author{Baofeng Wu}

\author{Dongdai Lin}

\address{State Key Laboratory of Information Security, Institute of Information Engineering, Chinese Academy of Sciences, Beijing 100093, China}

\begin{abstract}
We propose several techniques to construct complete permutation
polynomials of finite fields by virtue of complete permutations of
subfields. In some special cases, any complete permutation
polynomials over a finite field can be used to construct complete
permutations of certain extension fields with these techniques. The
results generalize some recent work of several authors.
\end{abstract}

\begin{keyword}
Permutation polynomial;  Complete; Subfield; Norm; Trace.\\
\medskip
\textit{MSC:~} 05A05 $\cdot$ 11T06 $\cdot$ 11T55
\end{keyword}

\end{frontmatter}

%-------------------------------------------------------------------------------------------------------------------

\section{Introduction}\label{secintro}

Let $\mathbb{F}_{q}$ be the finite field with $q$ elements where $q$
is a prime or a primer power. A polynomial $f(x)$ over
$\mathbb{F}_q[x]$ is called a permutation polynomial if the induced
map $\alpha\mapsto f(\alpha)$  from $\mathbb{F}_{q}$ to itself is
bijective \cite{lidl}. Note that it is only needed to study
permutation polynomials of normalized forms, i.e. permutation
polynomials which are monic and with no constant terms. Permutation
polynomials have important applications in combinatorics, coding and
cryptography. However, to construct large classes of them is far
from a simple matter. We refer to
\cite{akbary,charpin,xdhou,zbzha,zieve}, for example, for some
recent results on this topic.

A permutation polynomial $f(x)$ over $\mathbb{F}_{q}$ is further
called a complete one if $f(x)+x$ also plays as a permutation
polynomial. Complete permutation polynomials can be related to such
important combinatorial objects as orthogonal  latin squares.
However, it is of course more difficult to construct complete
permutation polynomials than constructing general permutation
polynomials. There were only a limited number of constructions known
before some great progresses were made recently.

In the recent work \cite{tu,xu,wu3,gfwu}, several new classes of
complete permutation polynomials were constructed. Some of these
classes are formed by monomials and the complete permutation
property of them can actually be characterized by a powerful lemma
given by  Zieve (see \cite[Lemma~2.1]{zieve4}); other classes are
formed by multinomials and the complete permutation property of them
can actually be characterized by the following general lemma.
\begin{lem}[See {\cite[Lemma~1.1]{akbary}}]\label{agwlem}
Let $A$, $S$ and $\bar{S}$ be finite sets with $\#S=\#\bar{S}$, and
let $f:A\rightarrow A$, $h:S\rightarrow\bar{S}$,
$\lambda:A\rightarrow S$ and $\bar{\lambda}:A\rightarrow \bar{S}$ be
maps such that $\bar{\lambda}\circ f=h\circ\lambda$. If both
$\lambda$ and $\bar{\lambda} $ are surjective, then the following
statements are
equivalent:\\
(1) $f$ is bijective (a permutation of $A$); and\\
(2) $h$ is bijective from $S$ to $\bar{S}$ and $f$ is injective on
$\lambda^{-1}(s)$ for each $s\in S$.
\end{lem}
\noindent In fact, a common feature of these complete permutation
polynomials is that they are related to permutations of certain
subsets of the finite fields, say subfields of the finite fields or
subgroups of multiplication groups of the finite fields. In
\cite{zieve3,zieve2} Zieve also generally studied how to induce
permutation polynomials over finite fields from subsets of them.

In the present paper, we generally study how to induce complete
permutation polynomials over finite fields from subsets of them. We
mainly focus on the case that the subsets are subfields of the
finite fields. More precisely, we propose several techniques to
construct complete permutation polynomials over finite fields by
virtue of complete permutation polynomials over certain subfields.
In some special cases, any complete permutations of the subfields
can be utilized. We also give some classes of complete permutation
polynomials constructed by these techniques which generalize some
recent work in \cite{tu,wu3,zieve2,gfwu}. The main tool that will be
used in our proofs is  Lemma \ref{agwlem} with $A=\mathbb{F}_{q^n}$
and $S=\bar{S}=\mathbb{F}_{q}$ where $n$ is a positive integer.

The rest of the paper is organized as follows. In Section
\ref{secnorm} and Section \ref{sectrace} we show how to construct
complete permutations of $\mathbb{F}_{q^n}$ from complete
permutations of $\mathbb{F}_{q}$ by fixing
$\lambda=\bar{\lambda}=\nor_{\mathbb{F}_{q^n}/\mathbb{F}_{q}}$ and
$\lambda=\bar{\lambda}=\tr_{\mathbb{F}_{q^n}/\mathbb{F}_{q}}$,
respectively, in Lemma \ref{agwlem}, where
$\nor_{\mathbb{F}_{q^n}/\mathbb{F}_{q}}$ and
$\tr_{\mathbb{F}_{q^n}/\mathbb{F}_{q}}$ represent the norm map and
trace map from $\mathbb{F}_{q^n}$ to $\mathbb{F}_{q}$, respectively.
Concluding remarks are given in Section \ref{seccon}.

%-------------------------------------------------------------------------------------------------------------------
\section{Complete permutation
polynomials induced from complete permutations of subfields and the
norm map}\label{secnorm}

In the sequel we fix an extension field $\mathbb{F}_{q^n}$ of the
finite field $\mathbb{F}_{q}$ and denote the norm and trace map from
$\mathbb{F}_{q^n}$ to $\mathbb{F}_{q}$ to be ``$\nor$" and ``$\tr$"
respectively for simplicity, i.e. $\nor(x)=x^{(q^n-1)/(q-1)}$ and
$\tr(x)=\sum_{i=0}^{n-1}x^{q^i}$ for any $x\in\mathbb{F}_{q^n}$.

Based on \cite[Lemma~2.1]{zieve4},  Zieve constructed many classes
of permutation polynomials over $\mathbb{F}_{q^n}$ by virtue of
permutation polynomials over $\mathbb{F}_{q}$ in \cite{zieve2}, some
of which cover certain known ones obtained by other authors via
complicated methods. For example, he obtained the following result.

\begin{lem}[See {\cite[Corollary~1.3]{zieve2}}]\label{zievePP}
Pick any $h\in\mathbb{F}_q[x]$, let $r$, $n$, $n'$ be positive
integers such that $nn'\equiv1\mod (q-1)$. Then
$x^rh\left(x^{(q^n-1)/(q-1)}\right)$ permutes $\mathbb{F}_{q^n}$ if
and
only if \\
(1) $(r,(q^n-1)/(q-1))=1$; and\\
(2) $g(x)=x^{rn'}h(x)$ permutes $\mathbb{F}_q$.
\end{lem}

In fact, item (2) in Lemma \ref{zievePP} is equivalent to  that
$\tilde{g}(x)=x^rh(x^n)$ permutes $\mathbb{F}_q$. From Lemma
\ref{zievePP}, we can directly obtain the following  construction of
complete permutation polynomials over finite fields.

\begin{thm}\label{normCPP}
Pick any $h\in\mathbb{F}_q[x]$, let  $n$  be a positive integers
such that $(n,q-1)=1$. Then $xh(\nor(x))$ is a complete permutation
polynomial of $\mathbb{F}_{q^n}$ if and only if $xh(x^n)$ is a
complete permutation polynomial of $\mathbb{F}_q$.
\end{thm}
\begin{pf}
It is a direct consequence of Lemma \ref{zievePP}   that
$xh\left(x^{(q^n-1)/(q-1)}\right)$ and
$x\left(h\left(x^{(q^n-1)/(q-1)}\right)+1\right)$ are permutations
of $\mathbb{F}_{q^n}$ if and only if $xh(x^n)$ and $x(h(x^n)+1)$ are
permutations of $\mathbb{F}_q$.\qedd
\end{pf}

Theorem \ref{normCPP} shows how to obtain complete permutation
polynomials of $\mathbb{F}_{q^n}$ from complete permutation
polynomials of $\mathbb{F}_{q}$ when $(n,q-1)=1$: we only need to
write a normalized complete permutation polynomial of
$\mathbb{F}_{q}$ into the form $xh(x^n)$ for certain
$h(x)\in\mathbb{F}_{q}[x]$ (note that every normalized polynomial,
can be written into this form since $(n,q-1)=1$), and then replace
$x^n$ by $\nor(x)=x^{(q^n-1)/(q-1)}$.

A special case of Theorem \ref{normCPP} was already obtained in
\cite{zieve2}, paying particular attention to finding complete
permutation monomials.

\begin{cor}\label{zieveCPP}
Pick $\alpha\in\mathbb{F}_q^*$ and let $n$ and $s$ be positive
integers with $(n,q-1)=1$. Then $f(x)=\alpha x^{1+s(q^n-1)/(q-1)}$
is a complete permutation polynomial of $\mathbb{F}_{q^n}$ if and
only if $\alpha x^{1+ns}$ is a complete permutation polynomial of
$\mathbb{F}_q$.
\end{cor}

\begin{rmk}
In fact the version of Corollary \ref{zieveCPP} in \cite{zieve2} is
that $f(x)$ is a complete permutation polynomial of
$\mathbb{F}_{q^n}$ if and only if $(1+s(q^n-1)/(q-1),q-1)=1$ and
$\alpha x^{1+ns}+x$ permutes $\mathbb{F}_q$. Obviously this is
equivalent to say that $\alpha x^{1+ns}$ is a complete permutation
polynomial of $\mathbb{F}_q$.
\end{rmk}

As an example, we give the following construction of complete
permutation polynomials which can be viewed as a generalization of
the construction in \cite[Corollary~3.2]{zieve2}. In fact,
\cite[Corollary~3.2]{zieve2} gives an answer to the open problem
posed in \cite{wu3}, thus our construction provides more examples
answering that open problem.

\begin{cor}\label{CPPeg}
Assume $q=r^t=2^{et}$. Let $k<t$ be a positive integer such that
$(k,t)\neq 1$ when $e=1$. Pick
$\alpha\in\mathbb{F}_q\backslash(\mathbb{F}_q)^{r^k-1}$. Then
$\alpha x^{1+(r^k-1)(q+1)q/2}$ is a complete permutation polynomial
of $\mathbb{F}_{q^2}$.
\end{cor}
\begin{pf}
Fix $n=2$ and $s=(r^k-1)q/2$ in Corollary \ref{zieveCPP}. It is
obvious that $(1+ns,q-1)=(r^k,q-1)=1$. Besides, $\alpha
x^{1+ns}+x=\alpha x^{r^k}+x$, which is a linearized permutation
polynomial of $\mathbb{F}_q$ since
$\alpha\in\mathbb{F}_q\backslash(\mathbb{F}_q)^{r^k-1}$.\qedd
\end{pf}

%-------------------------------------------------------------------------------------------------------------------
\section{Complete permutation
polynomials induced from complete permutations of subfields and the
trace map}\label{sectrace}

In fact, Theorem \ref{normCPP} can be derived form Lemma
\ref{agwlem} by fixing $A=\mathbb{F}_{q^n}$,
$S=\bar{S}=\mathbb{F}_{q}$ and $\lambda=\bar{\lambda}=\nor$. In this
section, we study how to obtain complete permutation polynomials
based on Lemma \ref{agwlem} by fixing $\lambda=\bar{\lambda}=\tr$.
As a matter of fact, many constructions of permutation polynomials
have already been obtained from Lemma \ref{agwlem} by fixing
$\lambda=\bar{\lambda}=\tr$ in \cite{yuan}. The following is an
example.

\begin{lem}[See {\cite[Corollory~5.2]{yuan}}]\label{yuanPP}
Pick $h(x)\in\mathbb{F}_q[x]$ with $h(0)\neq0$. Then $xh(\tr(x))$ is
a permutation polynomial of $\mathbb{F}_{q^n}$ if and only if
$xh(x)$ is a permutation polynomial of $\mathbb{F}_{q}$.
\end{lem}

Similar to Theorem \ref{normCPP}, the following construction of
complete permutation polynomials over finite fields can be obtained
directly from Lemma \ref{yuanPP}.

\begin{thm}\label{traceCPP1}
Pick $h(x)\in\mathbb{F}_q[x]$ with $h(0)\neq0,-1$. Then $xh(\tr(x))$
is a complete permutation polynomial of $\mathbb{F}_{q^n}$ if and
only if $xh(x)$ is a complete permutation polynomial of
$\mathbb{F}_{q}$.
\end{thm}
\begin{pf}
It is direct from Lemma \ref{yuanPP} that $xh(\tr(x))$ and
$x(h(\tr(x))+1)$ are permutations of $\mathbb{F}_{q^n}$ if and only
if $xh(x)$ and $x(h(x)+1)$ are permutations of $\mathbb{F}_{q}$
respectively.\qedd
\end{pf}

From Theorem \ref{traceCPP1} we know that by picking any normalized
complete permutation of $\mathbb{F}_{q}$ satisfying $h(0)\neq0,-1$
when written in the form $xh(x)$ for some $h(x)\in\mathbb{F}_q[x]$,
we can obtain a complete permutation of $\mathbb{F}_{q^n}$ by
replacing $x$ in the representation of $h(x)$ by $\tr(x)$. However,
when the condition $h(0)\neq0,-1$ is not satisfied, the method
cannot be applied directly. In the following, we propose another
construction of complete permutations of $\mathbb{F}_{q^n}$ based on
complete permutations of $\mathbb{F}_{q}$. We denote the kernel
space of the trace map from $\mathbb{F}_{q^n}$ to $\mathbb{F}_{q}$
to be $\ker(\tr)$ and assume $q=p^r$ where $p$ is a prime.

\begin{thm}\label{traceCPP2}
Pick $h(x)\in\mathbb{F}_{q}[x]$, let $L(x)$ be a $p$-polynomial over
$\mathbb{F}_{q}$, i.e. $L(x)$ is of the form
$\sum_{i=0}^{rn-1}a_ix^{p^i}$, $a_i\in\mathbb{F}_{q}^*$, $0\leq
i\leq r-1$, and denote $A(x)=L(x)/x$. Assume for some fixed
$a\in\mathbb{F}_{q}$, $L(x)-(h(b)/a+A(b))x$ and
$L(x)-[(h(b)+1)/a+A(b)]x$  can  induce  permutations of $\ker(\tr)$
for any $b\in\mathbb{F}_{q}$. Then the polynomial $xH(x)$ is a
complete permutation polynomial of $\mathbb{F}_{q^n}$ if and only if
$xh(x)$ is a complete permutation polynomial of $\mathbb{F}_{q}$
where
\[H(x)=h(\tr(x))+aA(\tr(x))-aA(x). \]
\end{thm}
\begin{pf}
It is obvious that
\begin{eqnarray*}
   \tr(x H(x))&=&\tr(x)h(\tr(x))+a\tr(x)A(\tr(x))-a\tr(xA(x))  \\
   &=&\tr(x)h(\tr(x))+aL(\tr(x))-a\tr(L(x))  \\
   &=&\tr(x)h(\tr(x))
\end{eqnarray*}
since $L(\tr(x))=\tr(L(x))$. Besides, for any $b\in\mathbb{F}_{q}$,
the limitation of $xH(x)$ on
$\tr^{-1}(b):=\{x\in\mathbb{F}_{q^n}\mid\tr(x)=b\}$ is
$x(h(b)+aA(b))-aL(x)=-a[L(x)-(h(b)/a+A(b))x]$, which can induce an
injective map on $\tr^{-1}(b)$ since $L(x)-(h(b)/a+A(b))x$ is
additive and can induce a permutation of $\ker(\tr)$. Applying Lemma
\ref{agwlem} by fixing $A=\mathbb{F}_{q^n}$,
$S=\bar{S}=\mathbb{F}_{q}$ and $\lambda=\bar{\lambda}=\tr$, we are
clear that $xH(x)$ is a permutation of $\mathbb{F}_{q^n}$ if and
only $xh(x)$ is a permutation of $\mathbb{F}_{q}$. By similar
arguments it follows that $x(H(x)+1)$ is a permutation of
$\mathbb{F}_{q^n}$ if and only $x(h(x)+1)$ is a permutation of
$\mathbb{F}_{q}$.\qedd
\end{pf}

Although the conditions in Theorem \ref{traceCPP2} are complicated,
explicit constructions belonging to Theorem \ref{traceCPP2} can be
obtained in some special cases. In fact, the main difficulty in
applying Theorem \ref{traceCPP2} is to construct permutation
$p$-polynomials of $\ker(\tr)$. A general discussion on constructing
such polynomials can be found in \cite{wuthesis}, which will not be
included in the present paper. We just focus on a special case in
the following, proposing a result which generalizes some recent work
in \cite{gfwu}.

\begin{lem}\label{binomialLPP}
Let $k$ be a positive integer with $(k,n)=1$ and $c\in\mathbb{F}_q$.
Then $x^{p^k}-cx$ permutes $\ker(\tr)$ in either
of the following cases:\\
(1)  $c^{(q-1)/(p^{(k,r)}-1)}=1$ and
$p\nmid n$; or\\
(2) $c^{n(q-1)/(p^{(k,r)}-1)}\neq1$.
\end{lem}
\begin{pf}
Since $x^{p^k}-cx\in\ker(\tr)$ for any $x\in\ker(\tr)$, we just need
to prove that $x^{p^k}-cx$ can induce an injective map on
$\ker(\tr)$ in either of the two cases.

In case (1), from $c^{(q-1)/(p^{(k,r)}-1)}=1$, i.e.
$c\in\left(\mathbb{F}_q^*\right)^{p^k-1}$, we know that the equation
$x^{p^k}-cx=0$ has nonzero solutions in $\mathbb{F}_q$. Besides,
from $(k,n)=1$ we have
\[(p^k-1,p^{nr}-1)=p^{(k,nr)}-1=p^{(k,r)}-1\mid(p^r-1),\]
thus all nonzero solutions of $x^{p^k}-cx=0$ in $\mathbb{F}_{q^n}$
are contained in $\mathbb{F}_q$. Therefore, from
$\mathbb{F}_q\cap\ker(\tr)=\{0\}$, which is implied by $p\nmid n$,
it follows that the system of equations $$\left\{\begin{array}{l}
x^{p^k}-cx=0\\\tr(x)=0
\end{array}\right.$$ has a unique solution $x=0$ in $\mathbb{F}_{q^n}$.
This is equivalent to say that $x^{p^k}-cx$ can induce an injective
map on $\ker(\tr)$;

In case (2), $x^{p^k}-cx$ can permute $\mathbb{F}_{q^n}$ since
$$c^{(q^n-1)/(q^n-1,p^k-1)}=c^{(q^n-1)/p^{(k,r)-1}}=c^{[(q^n-1)/(q-1)]\cdot[(q-1)/p^{(k,r)-1]}}=c^{n(q-1)/(p^{(k,r)}-1)}\neq1,$$ hence it can induce an injective
map on $\ker(\tr)$.\qedd
\end{pf}

\begin{rmk}
In fact, if $c^{(q-1)/(p^{(k,r)}-1)}=1$, then $p\nmid n$ is also a
necessary condition for $x^{p^k}-cx$ to permutes $\ker(\tr)$. This
is because if $p\mid n$, we have $\mathbb{F}_q\subseteq\ker(\tr)$.
However, $x^{p^k}-cx$ cannot induce a permutation of $\mathbb{F}_q$
due to $c^{(q-1)/(p^{(k,r)}-1)}=1$.
\end{rmk}

\begin{cor}\label{binomialLPP2}
Let $k$ be a positive integer with $(k,n)=1$. Assume $p\nmid n$ and
$(n,p^{(k,r)}-1)=1$. Then $x^{p^k}-cx$ permutes $\ker(\tr)$ for any
$c\in\mathbb{F}_q$.
\end{cor}
\begin{pf}
As $p\nmid n$, $x^{p^k}-cx$ can permutes $\ker(\tr)$ when
$c^{(q-1)/(p^{(k,r)}-1)}=1$ according to Lemma \ref{binomialLPP}
(1). Besides, obviously $x^{p^k}-cx$ can permutes $\ker(\tr)$ when
$c=0$. When $c\in\mathbb{F}_q^*$, $c^{(q-1)/(p^{(k,r)}-1)}\neq1$, we
claim that $c^{n(q-1)/(p^{(k,r)}-1)}\neq1$. This is because if
$c^{n(q-1)/(p^{(k,r)}-1)}=1$, we have
$c^{(n(q-1)/(p^{(k,r)}-1),q-1)}=1$. However, since
\[\left(n\frac{q-1}{p^{(k,r)}-1},q-1\right)=
\frac{q-1}{p^{(k,r)}-1}\left(n,p^{(k,r)}-1\right)=\frac{q-1}{p^{(k,r)}-1},\]
it follows that $c^{(q-1)/(p^{(k,r)}-1)}=1$, which leads to a
contradiction. Thus $x^{p^k}-cx$ can permutes $\ker(\tr)$ by Lemma
\ref{binomialLPP} (2).\qedd
\end{pf}

\begin{thm}\label{traceCPPeg}
Pick $h(x)\in\mathbb{F}_{q}[x]$, let $k$ be a positive integer with
$(k,n)=1$, and assume $p\nmid n$ and $(n,p^{(k,r)}-1)=1$. Then the
polynomial
\[x\left(h(\tr(x))+a\tr(x)^{p^k-1}-ax^{p^k-1}\right)\]
is a complete permutation polynomial of $\mathbb{F}_{q^n}$ if and
only if $xh(x)$ is a complete permutation polynomial of
$\mathbb{F}_{q}$ for any $a\in\mathbb{F}_{q}^*$.
\end{thm}
\begin{pf}
The result can be obtained directly by fixing $L(x)=x^{p^k}$ in
Theorem \ref{traceCPP2} and applying Corollary
\ref{binomialLPP2}.\qedd
\end{pf}

\begin{rmk}
Theorem \ref{traceCPPeg} generalizes \cite[Theorem~3]{chapuy} and
\cite[Theorem~5]{gfwu}.
\end{rmk}

By Theorem \ref{traceCPPeg}, we can use any complete permutation
polynomial of $\mathbb{F}_q$ to construct complete permutations of
certain extension fields.

%-------------------------------------------------------------------------------------------------------------------
\section{Concluding remarks}\label{seccon}
In this paper, we propose several general results on constructing
complete permutation polynomials of finite fields based on complete
permutations of subfields. Thanks to these results, many classes of
complete permutation polynomials over finite fields can be obtained
from known ones.

%-------------------------------------------------------------------------------------------------------------------

\end{document}